\documentclass{article}

\usepackage{amsmath}
\usepackage{amssymb}

\newtheorem{definition}{Definition}

\newtheorem{corollary}{Corollary}
\newtheorem{theorem}{Theorem}
\newtheorem{lemma}{Lemma}

\newtheorem{proposition}{Proposition}

\newcommand{\F}{{\cal F}}

\newcommand{\R}{{\bf R}}
\newcommand{\e}{\epsilon}
\newcommand{\N}{{\bf N}}
\newenvironment{proof}{{\par\noindent\bf Proof. }}{$\Box$}
\newenvironment{Proof}[1]{{\par\bf Proof of #1. }}{$\Box$}
\usepackage{latexsym}

\title{Measures of $\e$-complexity}
\author{V.Afraimovich and L.Glebsky\\ 
IICO-UASLP, A. Obregon 64
      San Lu\'\i s Potos\'\i, SLP 7820, Mexico}

\begin{document}
\maketitle
\abstract{We study some measures which are related to the notion of the 
$\epsilon$-complexity. We prove that measure of $\e$-complexity defined on
the base of the notion of $\e$-separability is equivalent to the dual measure that is
defined through $\e$-nets. 

}
\section{Introduction}

The problems under consideration in this article were originated 
in the process of study of complexity of behavior of orbits in dynamical systems.
While symbolic complexity (see, for instance \cite{F}) deals with symbolic systems and
topological complexity (\cite{BHM}) reflects pure topological features of dynamics, 
the $\e$-complexity depends essentially on a distance in the phase space (see definition bellow).
If one has a dynamical system generated by a continuous map $f:X\to X$ where $X$ is a metric space
with a distance $\rho$, one can introduce the sequence of distances
(\cite{B})
$$
\rho_n(x,y)=\max_{0\leq i\leq n-1} \rho (f^ix,f^iy),\ \ \ n\in \N,
$$
and study the $\e$-complexity with respect to the distance $\rho_n$ as a function of ``time'' $n$.
This function reflects the evolution of instability of orbits in time \cite{AZ}. 
But to study it in details, one
needs to know more about general properties of the $\e$-complexity of a metric space
(without dynamics).

The goal this article is to introduce and study quantities which contain an essential information 
about $\e$-complexity, the measures of $\e$-complexity in an ``abstract'' metric space.
The main results will be related to the $\e$-complexity defined on the base of the notion of 
$\e$-separability. The notion was used first by Kolmogorov and Tikhomirov \cite{KT} in their 
study of 
solutions of PDE and realization of random processes (Shannon suggested to pay attention to this 
notions in 1949, though). We will also study $\e$-complexities based on the notion 
of $\e$-nets.
We prove that measure of $\e$-complexity defined on
the base of the notion of $\e$-separability is equivalent to the dual measure that is
defined through $\e$-nets.

It appeared naturally that some results and ideas from discrete mathematics 
are worth to be exploited. We believe that we made the first step in this direction.

\section{Set-up and definitions}
\subsection{Separated sets and complexity}
Let $X,d$ be a compact metric space  with a distance $d$.
\begin{definition}
\begin{enumerate}
\item
Given $\epsilon>0$, a set $Y\subseteq X$ is $\epsilon$-separated iff
for any different $x,y\in Y$ one has $d(x,y)\geq\epsilon$. 
\item
The number 
$$
C_\epsilon (X,d)=
C_\epsilon:=\max\{|Y|,\;\; Y\mbox{is an $\epsilon$-separated set}\},
$$
where $|\cdot|$ denotes the cardinality of a set,
is called the $\epsilon$-complexity of $X$. 
\item An $\epsilon$-separated set $Y$ is optimal iff $|Y|=C_\epsilon$.
\end{enumerate}  
\end{definition} 
Let us show the following natural inequality.
\begin{proposition} \label{semiadditivity1}
Given $D_1,D_2\subseteq X$ and $\e>0$ one has
$$
C_\e(D_1\cup D_2)\leq C_\e(D_1)+C_\e (D_2).
$$
\end{proposition}
\begin{proof}
Let $Y\subseteq D_1\cup D_2$ be an optimal $\e$-separated set in $D_1\cup D_2$. Then $Y_i=Y\cap D_i$ is
an $\e$-separated set in $D_i$ and $|Y|\leq |Y_1|+|Y_2|\leq C_\e(D_1)+C_\e(D_2)$.
\end{proof}

{\bf Remark.}\\ 
Invariant sets in dynamical systems can be treated as results of inductive procedures.
For example, the dynamical system generated by the map $f:\R\to\R$, 
$$
f(x)=\left\{ \begin{array}{ll} 3x, & x\leq 1/2, \\
                              3x-3, & x>1/2,   \end{array} \right.  
$$ 
has an invariant set $K$ containing all orbits belonging to the interval $[0,1]$.
One can see that $K$ is the one-third Cantor set, so that 
$$
K=\bigcap_{n=1}^\infty \bigcup_{(i_0...i_{n-1})} \Delta_{i_0...i_{n-1}},
$$ 
where $i_j\in \{0,1\}$, $\Delta_{i_0...i_{n-1}}$ are intervals of the length $3^{-n}$
arising on the $n$-th step of construction of the Cantor set. Therefore, if
$\e\approx 3^{-n}$ then $C_\e\approx 2^n=\{$the number of different words of length $n$ in the 
full shift with $2$ symbols$\}=e^{hn}$, where $h=\ln 2$ is the topological entropy of the full 
shift.
Thus, 
$$
\frac{\ln C_\e}{-\ln \e}\approx \frac{\ln 2}{-\ln 1/3}=\frac{h}{\ln\lambda}=\dim_H K,
$$
where $\dim_H K$ is the Hausdorff dimension of $K$ and $\lambda=1/3$ is the contraction coefficient.
We obtained the familiar Furstenberg formula \cite{Fu}.

This example shows that if a subset of a metric space is the result of an inductive procedure governed
by a symbolic dynamical system then the $\e$-complexity contains, in fact, an important 
dynamical information. 
\subsection{$\e$-nets and complexity}
In this subsection we give a dual definition of complexity.
Given $x\in X$ let $O_\epsilon (x)=\{y\; :\; d(x,y)<\epsilon\}$, the ball of radius $\epsilon$
centered at $x$. Given $Y\subseteq X$ let $O_\epsilon(Y)=\bigcup\limits_{x\in Y} O_\epsilon(x)$. 
\begin{definition}
\begin{enumerate}
\item
Given $\epsilon>0$, a set $Y\subseteq X$ is an $\epsilon$-net iff
$O_\e(Y)=X$. 
\item
The number 
$$
R_\epsilon (X,d)=
R_\epsilon:=\min\{|Y|,\;\; Y\mbox{is an $\epsilon$-net}\},
$$
is called the dual $\epsilon$-complexity of $X$. 
\item An $\epsilon$-net  $Y$ is optimal iff $|Y|=R_\epsilon$.
\end{enumerate}
\end{definition}
The similar results to the one in Proposition~\ref{semiadditivity1} holds for dual complexities.
\begin{proposition}
Given $D_1,D_2\subseteq X$ and $\e>0$ one has
$$
R_\e(D_1\cup D_2)\leq R_\e(D_1)+R_\e (D_2).
$$
\end{proposition}
\begin{proof}
Let $Y_i\subseteq D_i$ be an optimal $\e$-net in $D_i$. Then $Y=Y_1\cup Y_2$ is
an $\e$-net in $D_1\cup D_2$ and $R_\e(D_1\cup D_2)\leq |Y|\leq |Y_1|+|Y_2|=R_\e(D_1)+R_\e(D_2)$.
\end{proof}

Any optimal $\e$-separated set is an $\e$ net, therefore $C_\e\geq R_\e$.
On the other hand the following statement holds.
\begin{proposition}
$R_{\e/2}\geq C_\e$
\end{proposition}
\begin{proof}
It follows directly from the definition that any pair of different points in an $\e$-separated
set $Z$ can not belong to a ball of radius $\e/2$. Thus we cannot cover $Z$ by less than 
$|Z|$ balls of radius $\e/2$. Assuming that $Z$ is optimal we obtain the inequality above.
\end{proof}

Let us introduce 
$$
b_\e=\sup_{x\in X} R_{\e/2}(O_\e(x)).
$$
Obviously, for any $D\subseteq X$ one has  $b_\e R_\e(D)\geq R_{\e/2}(D)$. 
It is not difficult to check that $b_\e\leq 2^d(2^d+1)$
for a subset of the Euclidean space $\R^d$.

\subsection{Ultrafilters}
Now we give some known results and definitions that can be found, for instance, in \cite{bu}.
\begin{definition}
A set $\F\subset 2^\N$ is called to be a filter over $\N$ iff it satisfies the following conditions:
\begin{itemize}
\item If $A\in \F$ and $B\in \F$, then $A\cap B\in\F$,
\item If $A\in \F$ and $A\subset B$ then $B\in \F$, 
\item $\emptyset\not\in\F$.
\end{itemize}
\end{definition} 
Let $a_n$ be a sequences of real numbers, $a$ is called to be a limit of $a_n$ with respect to
a filter $\F$, $a=\lim_\F a_n$, if for any $\e>0$ one has $\{n\ | \ |a_n-a|<\e\}\in\F$.
From the definition of a filter it follows that $\lim_\F a_n$ is unique, if exists.\\
{\bf Example} Let $\F_F=\{A\subseteq \N\ | \ \N\backslash A$ is finite $\}$. $\F_F$ is said to be a Frech\'et
filter. One can check that it is, indeed, a filter. A limit with respect to $\F_F$ coincides
with ordinary limit.
\begin{definition}
A filter $\F$ is called to be ultrafilter iff for any set $A\subseteq \N$ one has
$A\in \F$ or $\N\backslash A\in \F$. 
\end{definition}  
\begin{theorem}
A bounded sequences has a limit with respect to an ultrafilter. This limit is unique.
\end{theorem}
{\bf Example} For $i\in\N$ let  $\F_i=\{A\subseteq \N\ | \ i\in A\}$. It is an ultrafilter. 
Such an ultrafilter is called proper for $i$. 
One can check that 
$\lim_{\F_i} a_n=a_i$. So, limits with respect to a proper ultrafilter are not interesting.
\begin{proposition}
An ultrafilter $\F$ is proper (for some $i\in \N$) if and only if it contains a finite set.
\end{proposition}
This proposition implies that an ultrafilter is non-proper if and only if it is an extension
of the Frech\'et filter $\F_F$.
On the other hand, it follows  from the Zorn lemma that any filter can be extended to an ultrafilter.
\begin{proposition}
There is an ultrafilter $\F\supset \F_F$. Any such an ultrafilter is non-proper.
\end{proposition} 

\section{Measures of complexity} \label{complexity}

Our goal is to define a measure reflecting an asymptotic behavior of the 
$\epsilon$-complexity
as $\epsilon$ goes to $0$. For that we will use the technique of ultrafilters.

Given $\epsilon>0$, consider an optimal $\epsilon$-separated set
$A_\epsilon$. Introduce the following functional
$$
I_\epsilon(\phi)=\frac{1}{C_{\epsilon}}\sum_{x\in A_\epsilon}\phi(x)
$$
where $\phi:X\to R$ is a continuous function.
It is clear that $I_\epsilon$ is a positive bounded linear functional on $C(X)$.
Moreover, for any $\phi\in C(X)$ the family $I_\epsilon(\phi)$ is bounded. 
Fix a sequence $E=\{\e_n\}$, $\e_n\to 0$ as $n\to\infty$ and
an arbitrary non-proper ultrafilter $\F$. Consider
$$
I(\phi)=\lim_{\F}I_{\e_n}(\phi).
$$
$I$ is  a positive bounded linear functional on $C(X)$.
\begin{theorem} \label{independence_of_A}
The functional $I$ is independent of the choice
of an optimal sets $A_\e$.
\end{theorem}
\begin{proof}
 The proof is based on the  following proposition.
\begin{proposition}\label{bijection}
Let $A$ and $B$ be optimal $\epsilon$-separated
sets.
There exists a one-to-one map $\alpha: A\to B$ such that 
$d(x,\alpha(x))\leq \epsilon$ for any $x\in A$. 
\end{proposition}

Let $A_\e$ and $B_\e$ be optimal $\epsilon$-separated
sets, $\e\in E$. Let $\alpha_\e:A_\e\to B_\e$ be the map from Proposition~\ref{bijection}.
Then
$$
|\frac{1}{C_{\e}}\sum_{x\in A_\e}\phi(x)-
\frac{1}{C_{\e}}\sum_{x\in B_\e}\phi(x)|=
|\frac{1}{C_{\e}}\sum_{x\in A_\e}\left(\phi(x)-
\phi(\alpha_\e(x))\right)|\leq r_\phi (\epsilon)
$$
where $r_\phi(\epsilon)=\sup\{|\phi(x)-\phi(y)|\; :\; d(x,y)<\epsilon\}$,
the modulus of continuity of $\phi$. Since $X$ is a compact,
$r_\phi (\epsilon)\to 0$ as $\epsilon \to 0$. It implies the desired result due to the choice
of the ultrafilter $\F$.
So, we need only to prove Proposition~\ref{bijection}; it will be done below.
\end{proof}

In the proof of Proposition~\ref{bijection} we will need the Marriage Lemma of P. Hall,
see for instance \cite{Ryser}.
\begin{lemma} \label{mariage}
For an indexed collections of finite sets $F_1, F_2,\ldots, F_k$ the following conditions
are equivalent:
\begin{itemize}
\item there exists an injective function $\alpha : \{1,2,...,k\}\to\bigcup\limits_{i=1}^k F_i$
such that $\alpha(i)\in F_i$;
\item For all $S\subseteq \{1,2,\ldots,k\}$ one has $|\bigcup\limits_{i\in S} F_i|\geq |S|$.
\end{itemize}
\end{lemma}

Recall that $O_\epsilon (x)=\{y\; :\; d(x,y)<\epsilon\}$, the ball of radius $\epsilon$
centered at $x$. Given $Y\subseteq X$ let $O_\epsilon(Y)=\bigcup\limits_{x\in Y} O_\epsilon(x)$. 
\begin{Proof}{Proposition~\ref{bijection}}
For any $x\in A$ let $B_x=O_\epsilon (x)\cap B$. If we show that for any $S\subseteq A$ the
following inequality holds
\begin{equation} \label{Hall_condition}
|\bigcup_{x\in S} B_x|\geq |S|,
\end{equation}
then the proposition follows from Lemma~\ref{mariage} due to $|A|=|B|=C_\epsilon$.
To prove inequalities~(\ref{Hall_condition}), suppose that 
$|\bigcup\limits_{x\in S} B_x|=|O_\epsilon(S)\cap B|< |S|$ for some $S\subseteq A$. Then
$$
|S\cup(B\setminus (O_\epsilon(S)\cap B)|=|S|+(|B|-|O_\epsilon(S)\cap B|)>|B|=C_\epsilon,
$$
on the other hand, the set $S\cup(B\setminus (O_\epsilon(S)\cap B)$ is $\epsilon$-separated.
We have a contradiction with optimality of $B$.
\end{Proof}

So, we have defined a functional $I$ which may depend on the choice of the sequence $E$ and the ultrafiter $\F$ 
only.
Sometimes we will write $I_{E,\F}$ to emphasize this dependence. It is well known, that $I_{E,\F}$
generate unique regular Borel measure $\mu_{E,\F}$ on $X$ such that $\mu_{E,\F}(X)=1$. 
\begin{definition}
The measures  $\mu_{E,\F}(X)$ will be called  measures of complexity.
\end{definition}
We are going to show examples of $(X,d)$ when $\mu_{E,\F}=\mu$ is independent on $E$, $\F$ and
when $\mu_{E,\F}$ depends on $E, \F$. In the first case
$$
I_{E,\F}(\phi)=I(\phi)=\lim_{\e\to 0}I_\e(\phi).
$$

Of course, it is difficult to find optimal sets and construct directly measures of complexity in real situations.
Nevertheless, it is possible to work with them by using some of their intrinsic properties. Let us show now that 
measures of complexity are invariant with respect to local isometries. 
\begin{definition}
A homeomorphism $\tau:X\to X$ is called to be $\e$-isometry iff
$d(x,y)=d(\tau(x),\tau(y))$ for all $x,y\in X$, $d(x,y)\leq\e$ 
A homeomorphism $\tau:X\to X$ is called to be local isometry iff it is
$\e$ isometry for some $\e>0$. 
\end{definition}        
 It is clear that an isometry is a local isometry.
\begin{proposition} \label{Prop_group}
Local isometries with composition form a group. 
\end{proposition}
\begin{proof}
It is easy to check that the
composition of two $\e$-isometries is an $\e$-isometry.
Let $\tau$ be an $\e$-isometry. Then $\tau^{-1}$ is uniformly 
continuous and there exists $\e'>0$ such that if $d(x,y)\leq\e'$,then
$d(\tau^{-1}(x),\tau^{-1}(y))\leq\e$. Consequently, if $d(x,y)\leq\e'$
then $d(x,y)=d(\tau^{-1}(x),\tau^{-1}(y))$, so, $\tau^{-1}$ is an
$\e'$-isometry.
\end{proof}

We do not know if $\e$-isometries form a group.

\begin{proposition}\label{Prop_invers_separation}
Let $\tau$ be an $\e_0$-isometry and $A$ be an $\e$-separated set, $\e\leq \e_0$.
Then $\tau^{-1}(A)$ is also $\e$-separated.
\end{proposition}
\begin{proof}
Assume, on the contrary, that $\tau^{-1}(A)$ is not $\e$-separated, i.e.,
there are different $x,y\in\tau^{-1}(A)$ with $d(x,y)<\e\leq \e_0$. Then
$d(\tau(x),\tau(y))=d(x,y)<\e$, so $A$ cannot be $\e$-separated.  
\end{proof}
\begin{theorem} \label{Th_invariance1}
Let $\tau$ be a local isometry. Then $\mu_{E,\F}$ is invariant,
i.e. $\mu_{E,\F}(A)=\mu_{E,\F}(\tau^{-1}(A))$ for all measurable $A$.
\end{theorem}
\begin{proof}
It is enough to show that for all $\phi\in C(X)$
\begin{equation} \label{eq_F}
I_\F(\phi\circ\tau)=I_\F(\phi).
\end{equation}
There exists $\e_0>0$, such that $\tau$ is an $\epsilon_0$-isometry. Let $A_\e$ be an optimal
$\e$-separated set, $\e\leq \e_0$. It follows from Proposition~\ref{Prop_invers_separation} 
that $\tau^{-1}(A_\e)$ is an optimal $\e$-separated set.
It implies the validity of Equation~(\ref{eq_F}). Indeed, 
$$
I_\epsilon(\phi)=\frac{1}{C_{\epsilon}}\sum_{x\in A_\epsilon}\phi(x),\ \ 
I_\epsilon(\phi\circ\tau)=\frac{1}{C_{\epsilon}}\sum_{x\in \tau^{-1}(A_\epsilon)}\phi(x),
$$ 
and the result follows from Theorem~\ref{independence_of_A}.
\end{proof}

\begin{corollary}
Let a continuous group  operation $*$ be defind on $X$ such that right shifts $r_g(x)=g*x$ 
(left shifts $l_g(x)=x*g$) are local isometries for all $g\in X$. Then $\mu_{E,\F}$ is
the normalized Haar measure on $(X,*)$. In particular, $\mu_{E,\F}$ does not depend on 
$E,\F$.   
\end{corollary}

{\bf Example 1.}
Let $X=\Omega_p$, the full shift with $p$ symbols, i.e. $\Omega_p=\{0,1,...,p-1\}^{Z^+}$  
with the distance
$$
d_q(x,y)=\sum_{i=0}^\infty \frac{|x_i-y_i|}{q^i},\; q>1.  
$$
$\Omega_p$ can be equipped by the group operation $\oplus$ as follows:
$$
(x\oplus y)_i=x_i+y_i, \mod p
$$
It is clear that $(\Omega_p,\oplus)$ is a continuous group. Moreover, the right translation
by any element is an isometry. Therefore $\mu_{E,\F}=\mu$ coincides with the Haar measure which, 
in fact, is the $(1/p,...,1/p)$-Bernoulli measure.

{\bf Example 2.}
Let $X=\Omega_{M}$ be a topological Markov chain, defined by a finite matrix 
$M:\{0,1,...,p-1\}^2\to\{0,1\}$,
i.e. $\Omega_{M}=\{<x_0,x_1,...>\;|\; x_i\in \{1,2,...,p-1\}\;\mbox{and}\; M(x_i,x_{i+1})=1\}$.
Metric $d$ is the same as in  Example~1.

Cylinder $[a_0,a_1,...,a_{n-1}]$ of the length $n$ is the set of all $x\in\Omega_M$,
such that $x_i=a_i$ for $i=0,1,...,n-1$. 
A word $<a_0,a_1,...,a_{n-1}>$ is 
admissible iff  $[a_0,a_2,...,a_{n-1}]\neq\emptyset$. Let $W_n$ be the set of all admissible
words of the length $n$ and $\alpha$ be a permutation of $W_n$ such that 
$(\alpha(w))_{n-1}=w_{n-1}$ for every $w\in W_n$ (admissible permutation). 
Given such an $\alpha$ define 
$g_\alpha:X\to X$ as follows
$$
g_\alpha(x)=(\alpha(x_0,x_1,...,x_{n-1}),x_n,x_{n+1},...).
$$  
It is simple to see that $g_\alpha$ is a local isometry. It implies that
$\mu_{E,\F}([a_0,a_1,...,a_{n-1}])=\mu_{E,\F}([b_0,b_1,...,b_{n-1}])$ if
$[a_0,a_1,...,a_{n-1}]\neq\emptyset$, $[b_0,b_1,...,b_{n-1}]\neq\emptyset$ and 
$a_{n-1}=b_{n-1}$. Indeed, under these assumptions there exists an admissible permutation 
$\alpha:W_n\to W_n$ such that $\alpha(a_0,a_1,...,a_{n-1})=b_0,b_1,...,b_{n-1}$. 
So, the measure $\mu_{E,\F}$ of a nonempty cylinder $[a_0,a_1,...,a_{n-1}]$ depends only on 
$a_{n-1}$ and $n$. Let $v_i(n)=\mu_{E,\F}([a_0,a_1,...,a_{n-2},i])$ for an admissible
$<a_0,a_1,...,a_{n-2},i>$ ($v_i(n)=0$ if there is no admissible words of length $n$ ending by $i$).
It is simple to check that 
$$
v_i(n)=\sum_{j,M(i,j)=1}v_j(n+1)
$$
This relation can be rewritten in the matrix form
$$
v(n)=Mv(n+1),
$$
where $v(n)=(v_0(n),v_1(n),...,v_{p-1}(n))^T$ is a column vector.  
If $M$ is a primitive matrix ($M^p>0$ for some $p$) then this
equation uniquely defines the measure $\mu_{E,\F}$, which 
in this case 
turns out to be independent of $E,\F$. Indeed, by Perron Theorem
matrix $M$ has unique positive eigenvector $e$ with eigenvalue
$\lambda>0$ (in our case, in fact, $\lambda>1$). Let $P$ be the 
set of all  lines in $\R^p$,
generated by non-negative vectors. From the proof of Perron Theorem (see, for example, \cite{Katok})  
$$
\bigcap_{n\in\N}M^n(P)= \{l_e\},
$$
where $l_e$ is a line, generated by $e$.
Since $v(n)>0$ and $v(k)=M^nv(n+k)$, one has 
$l_ {v(k)}\in M^n(P)$ for any $n$. Hence, $v(k)=c_ke$. So, 
$v(n)=\lambda^{-n}c_0e$. We have proved the following
\begin{proposition}
Let $M$ be a primitive matrix and 
$C\subset \Omega_M$ is an admissible cylinder of length $n$, ending by
$i$. Then 
$\mu_{E,\F}(C)=\lambda^{-n}e_i$, where $(e_0,e_1,...,e_{p-1})$ is the
positive eigenvector of $M$, with $e_0+e_1+...e_{p-1}=1$.
\end{proposition}

{\bf Example 3.} Here we construct an example where $\mu_{E,\F}$ is not unique.
Let $X=\Omega_{0,1}\cup\Omega_{2,3}$, where $\Omega_{i,j}$ is the 
Bernoulli shift of symbols $i,j$. We are going to introduce
a metric $d$ on $X$ such that $\mu_{E,\F}$ depends on $E.\F$.

Let us define $d$. For $x\in\Omega_{0,1}$ and $y\in\Omega_{2,3}$
let $d(x,y)=1$. For $x,y\in\Omega_{0,1}$, $x_n\neq y_n$ and $x_i=y_i$ 
for $i<n$, let $d(x,y)=a_n$. 
For $x,y\in\Omega_{2,3}$, $x_n\neq y_n$ and $x_i=y_i$ 
for $i<n$, let $d(x,y)=b_n$. 
Suppose, $1\geq a_0\geq a_1\geq...\geq a_n\to 0$ and 
$1\geq b_0\geq b_1\geq...\geq b_n\to 0$. 
Straightforward calculations show that $d$ is a metric (even an ultrametric) defining
the Markov topology on $X$. 

\begin{proposition}
If $a_{r-1}\geq\e> a_r$ and 
$b_{m-1}\geq\e > b_m$ then $C_\e(\Omega_{0,1})=2^r$ and 
$C_\e(\Omega_{2,3})=2^m$, the cardinality of an optimal $\e$-separated set
on $\Omega_{0,1}$ and $\Omega_{2,3}$, correspondingly. 
\end{proposition}
\begin{proof}
Indeed, if, say, $x,y\in \Omega_{0,1}$
are in the same cylinder of length $r$, then $d(x,y)\leq a_r<\e$. So, an $\e$-separated set
does not contain different points of the same cylinder of length $r$. On the other hand, 
if $x,y\in \Omega_{0,1}$
are in different cylinders of length $r$, then $d(x,y)\geq a_{r-1}\geq\e$.
\end{proof}  

Take $\e_n=1/2n$ and $\e'_n=1/(2n+1)$. The idea is to choose $a_n$ and 
$b_n$ such that 
\begin{equation}\label{eq_limits}
\frac{C_{\e_n}(\Omega_{2,3})}{C_{\e_n}(\Omega_{0,1})}\to 0\
\mbox{and} \ \frac{C_{\e'_n}(\Omega_{0,1})}{C_{\e'_n}(\Omega_{2,3})}\to 0,
\end{equation}
as $n\to\infty$. In particular, we can take $b_0=1$,
$a_{(n-1)(2(n-1)+1)}= a_{(n-1)(2(n-1)+1)+1}=...=a_{n(2n+1)-1}=1/2n$ and 
$b_{n(2n-1)}=b_{n(2n-1)+1}=...=b_{(n+1)(2n+1)-1}=1/(2n+1)$, where $n=1,2....$.
Now one can check that 
$$
a_{n(2n+1)-1}=\frac{1}{2n}=\e_n> a_{n(2n+1)}\ \ \mbox{and}\ \
b_{n(2n-1)-1}>\frac{1}{2n}=\e_n> b_{n(2n-1)}.
$$
Because of the proposition  $C_{\e_n}(\Omega_{0,1})=2^{n(2n+1)}$, 
$C_{\e_n}(\Omega_{2,3})=2^{n(2n-1)}$ and the first limit in (\ref{eq_limits}) occurs.

On the other hand
$$
a_{n(2n+1)-1}>\frac{1}{2n+1}=\e'_n> a_{n(2n+1)}\ \
b_{(n+1)(2n+1)-1}=\frac{1}{2n+1}=\e'_n> b_{(n+1)(2n+1)}.
$$
So, $C_{\e'_n}(\Omega_{0,1})=2^{n(2n+1)}$ and $C_{\e'_n}(\Omega_{2,3})=2^{(n+1)(2n+1)}$;
the second limit of (\ref{eq_limits}) is valid. Now,
for $E=\{1/2n\;|\;n\in\N\}$ one has
$\mu_{E,\F}(\Omega_{2,3})=0$ and $\mu_{E,\F} |_{\Omega_{0,1}}$ is 
the $(1/2,1/2)$-Bernoulli measure, independently of $\F$.  
For $E'=\{1/(2n+1)\;|\;n\in\N\}$ one has
$\mu_{E',\F}(\Omega_{0,1})=0$ and $\mu_{E',\F} |_{\Omega_{2,3}}$ is 
the $(1/2,1/2)$-Bernoulli measure, independently of $\F$.
For $\tilde E=\{1/n\; |\; n\in\N\}=E\cup E'$ the measure $\mu_{\tilde E,\F}$ will
depend on $\F$.
\section{Measures of dual complexity}
To define measures of dual complexity we proceed in the same way as in Section~\ref{complexity},
just replacing $\e$-separated sets by $\e$-nets. 

Given $\epsilon>0$, consider an optimal $\epsilon$-net
$A_\epsilon$. Introduce the following functional
$$
\tilde I_\epsilon(\phi)=\frac{1}{R_{\epsilon}}\sum_{x\in A_\epsilon}\phi(x).
$$
Consider
$$
\tilde I(\phi)=\lim_{\F}\tilde I_{\e_n}(\phi).
$$
\begin{theorem} \label{independence_of_B}
The functional $\tilde I$ is independent of the choice
of an optimal $\e$-nets $A_\e$.
\end{theorem}
\begin{proof}
The proof is similar to the one of Theorem~\ref{independence_of_A},
just instead of Proposition~\ref{bijection} one should use Proposition~\ref{bijection1},
formulated below.
\end{proof}
\begin{proposition}\label{bijection1}
Let $A$ be an optimal $\epsilon$-net and $B$ be an $\e$-net.
There exists an injective map $\alpha: A\to B$ such that 
$d(x,\alpha(x))\leq 2\epsilon$ for any $x\in A$. 
\end{proposition}
\begin{proof}
Again we will use Marriage Lemma (Lemma~\ref{mariage}). For $x\in A$ let
$$
B_x=\{y\in B\; |\; O_\e(y)\cap O_\e(x)\neq\emptyset\}
\subseteq O_{2\e}(x)\cap B
$$
For $S\subseteq A$ let 
$$
B_S=\bigcup_{x\in S} B_x.
$$
As in the proof of Proposition~\ref{bijection} it is 
enough to show that for any $S\subseteq A$ one has 
$$
|B_S|\geq |S|. \eqno(*)
$$
First of all, $O_\e(x)\subseteq O_\e(B_x)$, $x\in A$. Indeed, due to
$O_\e(B)=X$ we have
$$
O_\e(x)=O_\e(x)\cap O_\e(B)=O_\e(x)\cap O_e(B_x).
$$ 
So, $O_\e(S)\subseteq O_\e(B_S)$. 

Now, suppose that $|B_S|< |S|$  in contradiction 
to (*). Then 
$$
|A\backslash S\cup B_S|<|A|.
$$ 
Moreover, $O_\e(A\backslash S)\supseteq O_\e(A)\backslash O_\e(S)$. Indeed,
if $z\in O_\e(A))$ and $z\not\in O_\e(S)$ then there exists $a\in A$ such that $d(a,x)<\e$;
$a$ cannot belong to $S$ because $z\not\in O_\e(S)$. Hence, $a\in A\backslash S$, and 
$z\in O_\e(A\backslash S)$. Thus
$$
O_\e(A\backslash S\cup B_S)=O_\e(A\backslash S)
\cup O_\e(B_S)\supseteq 
O_\e(A)\backslash O_\e(S)\cup O_\e(B_S)=X\backslash O_\e(S)\cup O_\e(B_S) =X,
$$
the contradiction with minimality of $A$.
\end{proof}
\begin{definition}
The measures  $\nu_{E,\F}(X)$ corresponding to  $\tilde I_{E,\F}$ will be called dual  
measures of complexity.
\end{definition}
\begin{proposition}\label{Prop_separation_net}
Let $\tau$ be an $\e_0$-isometry and $A$ be an $\e$-net, $\e\leq \e_0$.
Then $\tau(A)$ is also an $\e$-net.
\end{proposition}
\begin{proof}
Given $x\in X$ we have to prove that $x\in O_\e(\tau(A))$.
Due to surjectivity of $\tau$ there exists $y\in X$, $x=\tau(y)$.
There exists $a\in A$ such that $y\in O_\e(a)$. By the definition of $\e$-isometry
$x=\tau (y)\in O_\e(\tau(a))$.  
\end{proof}

Using Proposition~\ref{Prop_separation_net}, Proposition~\ref{Prop_group}
one can prove the following analogue of Theorem~\ref{Th_invariance1}.
\begin{theorem} \label{Th_dual_invariant}
Let $\tau$ be a local isometry. Then $\nu_{E,\F}$ is invariant,
i.e. $\nu_{E,\F}(A)=\nu_{E,\F}(\tau^{-1}(A))$ for all measurable $A$.
\end{theorem}

We don't know if $\mu_{E,\F}$ and $\nu_{E,\F}$ can be different,
but we can prove the following theorem.
\begin{theorem}\label{Th_relation}
If there exists $k\in\N$ such that for any $x\in X$ and any small enough $\e>0$
one has $C_\e(O_\e(x))\leq k$, then $\mu$ and $\nu$ are equivalent and, moreover,
$$
\frac{1}{k}\nu_{E,\F}(A)\leq\mu_{E,\F}(A)\leq k\nu_{E,\F}(A)
$$

for any Borel set  $A\subseteq X$. 
\end{theorem} 
It easily implies
\begin{corollary}
If $d$ is a ultrametric (i.e.  $d(x,z)\leq\max\{d(x,y), d(y,z)\}$
for any $x,y,z\in X$) then $\nu_{E,\F}=\mu_{E,\F}$.
\end{corollary}
\begin{proof}
The result follows from the fact that $C_\e(O_\e(x))=1$ for any $x\in X$ and any $\e>0$, so
$k=1$ in the conditions of the theorem.
Indeed, for any $z,y\in O_\e(x)$ one has $d(z,y)\leq\max \{d(z,x),d(x,y)\}<\e$
\end{proof}

So, for Example~3 of Section~\ref{complexity} one has $\nu_{E,\F}=\mu_{E,\F}$.
The measures in Examples~1,2 are also coinside because of Theorem~\ref{Th_dual_invariant}.

In our proof of Theorem~\ref{Th_relation} we will use the following proposition.
\begin{proposition}\label{sets}
Let $A$ be an optimal $\e$-net and $B$ be an optimal $\e$-separated set.
Then there exists a collection $\{K_x\}$ of subsets of $B$, 
indexed by elements of $A$, with the following properties:
\begin{itemize}
\item $K_x\subseteq O_\e(x)\cap B$ for any $x\in A$;
\item $K_x\neq\emptyset$ for any $x\in A$;
\item $K_x\cap K_y=\emptyset$ for any different $x,y\in A$;
\item $\bigcup_{x\in A}K_x=B$.
\end{itemize} 
\end{proposition}
\begin{proof}
Since $B$ is an $\e$-net,
it follows from Proposition~\ref{bijection1} that there exists an injective map 
$\alpha: A\to B$. So, we can put $\alpha(x)$ to $K_x$ and distribute
the points $B\backslash \alpha(A)$ among $K_x$ so that $K_x$ satisfy the
properties claimed. (For example, we can order $A$ and put 
$b\in B\backslash \alpha(A)$ into $K_x$ with the smallest $x\in A$ such that $b\in O_\e(x)$).
\end{proof} 
\begin{Proof}{Theorem~\ref{Th_relation}}
It is enough to show that for non-negative continuous $\phi$
\begin{equation} \label{func_ineq}
\frac{1}{k}\tilde I_\e(\phi)-\delta_\e(\phi)\leq
I_\e(\phi) \leq k\tilde I_\e(\phi)+\delta_\e(\phi),
\end{equation}
where $\delta_\e(\phi)$ is the modulus of continuity of $\phi$.
Let $A$ be an optimal $\e$-net and $B$ be an optimal $\e$-separated set.
Let $K_x$ be the sets of Proposition~\ref{sets}. 
From the conditions of the theorem it follows that $|K_x|\leq k$ and 
$R_\e\leq C_\e\leq kR_\e$.
Then
$$
\sum_{y\in B}\phi(y)=\sum_{x\in A}\sum_{y\in K_x}\phi(y)\leq 
\sum_{x\in A} |K_x|(\phi(x)+\delta_\e(\phi))\leq k \sum_{x\in A} \phi(x)+ C_\e\delta_\e(\phi).
$$
Thus
$$
\frac{1}{C_\e}\sum_{y\in B}\phi(y)\leq 
\frac{k}{C_\e}\sum_{x\in A}\phi(x) +\delta_\e(\phi)\leq 
\frac{k}{R_\e}\sum_{x\in A}\phi(x) +\delta_\e(\phi),
$$
that proves the right inequality in (\ref{func_ineq}).  Similarly,
$$
\sum_{y\in B}\phi(y)=\sum_{x\in A}\sum_{y\in K_x}\phi(y)\geq 
\sum_{x\in A}|K_x| (\phi(x)-\delta_\e(\phi))\geq \sum_{x\in A} \phi(x)- C_\e\delta_\e(\phi).
$$
Thus
$$
\frac{1}{C_\e}\sum_{y\in B}\phi(y)\geq
\frac{1}{C_\e}\sum_{x\in A}\phi(x) -\delta_\e(\phi)\geq 
\frac{1}{kR_\e}\sum_{x\in A}\phi(x) -\delta_\e(\phi),
$$
\end{Proof}



\underline{\bf Acknowledgments}. The authors would like to thank P. Collet, E. Ugalde and G.M. Zaslavsky 
for useful discussions. V.A. was partially supported by CONACyT grant 485100-5-36445-E.
L.G. was partially supported by CONACyT-NSF grant E120.0547 and PROMEP grant PTC-62.

\end{document}